\documentclass[final,3p,authoryear]{elsarticle}

\usepackage{amssymb}

\usepackage[T1]{fontenc}
\usepackage[latin1]{inputenc}
\usepackage{amsmath}
\usepackage{mathrsfs}
\usepackage[colorlinks]{hyperref}

\newtheorem{theorem}{Theorem}
\newtheorem{lemma}{Lemma}

\newtheorem{proposition}{Proposition}
\newtheorem{remark}{Remark}

\newcommand{\al}{\alpha}
\newcommand{\be}{\begin{equation}}
\newcommand{\ee}{\end{equation}}
\newcommand{\half}{\frac{1}{2}}
\newcommand{\p}{\partial}

\begin{document}

\begin{frontmatter}

\title{A bending theory of thermoelastic diffusion plates based on Green-Naghdi theory}

\author[tun]{Moncef Aouadi\corref{cor1}}
\ead{moncefaouadi00@gmail.com}

\author[sal]{Francesca Passarella}
\ead{fpassarella@unisa.it}

\author[sal]{Vincenzo Tibullo}
\ead{vtibullo@unisa.it}

\cortext[cor1]{Corresponding author}

\address[tun]{Ecole Nationale d'Ingénieurs de Bizerte, Université de Carthage\\
BP66, Campus Universitaire, Menzel Abderrahman 7035, Tunisia}

\address[sal]{Università degli Studi di Salerno, Dipartimento di Matematica\\
Via Giovanni Paolo II 132, 84084 Fisciano (SA), Italy}

\begin{abstract}
This article is concerned with bending plate
 theory for thermoelastic diffusion materials under Green-Naghdi theory.
 First, we present the basic equations which characterize the bending
of thin thermoelastic diffusion plates for type II and III models.
The theory allows for the effect of transverse shear deformation
without any shear correction factor, and permits the propagation of
waves at a finite speed without energy dissipation for type II model
and with energy dissipation for type III model. By the semigroup
theory of linear operators, we prove the well-posedness of the both
models and the asymptotic behavior of the solutions of type III
model. For unbounded plate of type III model, we prove that a
measure associated with the thermodynamic process decays faster than
an exponential of a polynomial of second degree.
 Finally,  we investigate the
impossibility of the localization in time of  solutions. The main
idea to prove this result is to show the uniqueness of solutions for
the backward in-time problem.

Link to publisher: \url{https://doi.org/10.1016/j.euromechsol.2017.03.001}
\end{abstract}


\end{frontmatter}


\section{ Introduction}
Elastic plates play an important role in mechanical structures since
they can support loads far in excess of their own weight. In
addition, due to their geometric characteristics, thin plates can be
studied mathematically by means of two-dimensional models instead of
the full and much more complex equations of three-dimensional
elasticity.
 The first plate-bending model was
proposed by \cite{Kirchhoff}. Making a number of
simplifying hypotheses, he arrived at the conclusion that, in terms
of Cartesian coordinates $(x_1, x_2, x_3)$ with $(x_1, x_2)$ in the
middle plane of the plate, the displacement field should be of the
form $(x_3u_1, x_3u_2,u_3)$, where the functions $u_i =u_i(x_1,
x_2),\ i =1, 2, 3$, satisfy $u_1=-u_{3,1},\ u_2=-u_{3,2},\
\Delta\Delta u_3= q /D$, $q$  being the total load, $D$ the modulus
of rigidity of the material and $\Delta$ the two-dimensional
Laplacian. It is clear that the nonhomogeneous biharmonic equation
for $u_3$ cannot take more than two boundary conditions.

 Since there are numerous cases where the transverse
shear force is not negligible and each of the three moments must be
given on the contour, the need arose for more refined models with a
more sophisticated mathematical content. One such model was proposed
by \cite{Reissner3}, who started with the stress tensor,
postulating a linear dependence on $x_3$ for the components
$t_{\alpha\beta},\ \alpha,\beta\in \{1, 2\}$ and a certain type of
parabolic dependence on $x_3$ for the components $t_{\alpha 3},\
\alpha =1, 2$. This model accepts three independent boundary
conditions, but does not yield the explicit expressions of the
displacements. Another model, proposed by \cite{Mindlin}, is
based on Kirchhoff's kinematic assumption on the displacement field
but has no tie between its first two components and the third one,
as above; however, it makes use of a correction factor in the
constitutive equations, which interferes with its mathematical
rigor. Again, this model accepts three boundary conditions. It
should be pointed out that all Reissner and Mindlin type models also
account for the transverse shear force in the plate.

On the other hand, the temperature plays a significant role in the
process of bending thermoelastic plates. The deformation of such
plates is of interest in a wide variety of practical problems, from
microchip production to aerospace industry. Relevant theoretical
developments on the subject were made by 
\cite{Green2, Green3}. They developed three models for generalized
thermoelasticity of homogeneous isotropic materials which are
labeled as model I, II and III. This theory is developed in a
rational way to produce a fully consistent theory that is capable of
incorporating thermal pulse transmission in a very logical manner.
When the theory of type I is linearized, the parabolic equation of
the heat conduction arises. Type II theory predicts a finite speed
for heat propagation and involves no energy dissipation, now
referred to as thermoelasticity without energy dissipation. Type III
theory permits propagation of thermal signals at both finite and
infinite speeds.

We may think that the classical theory of thermoelastic bending
plates is a good model to explain the thermal conduction in
different structures. However, the research conducted in the
development of high technologies after the second world war,
confirmed that the field of diffusion in solids cannot be ignored.
Thus, the obvious question is, what happens when the diffusion
effect is considered with the thermal effect in the theory of
bending elastic plate. Diffusion can be defined as the random walk
of a set of particles from regions of high concentration to regions
of lower concentration. Thermodiffusion in an elastic solid is due
to coupling of the fields of strain, temperature and mass diffusion.
The processes of heat and mass diffusion play important roles in
many engineering applications, such as satellites problems,
returning space vehicles and aircraft landing on water or land.
There is now a great deal of interest in the process of diffusion in
the manufacturing of integrated circuits, integrated resistors,
semiconductor substrates and MOS transistors. Oil companies are also
interested to this phenomenon to improve the conditions of oil
extractions.

Recently, \cite{Aouadi} used the results of Green and
Naghdi on thermo-mechanics of continua to derive a nonlinear theory
of thermoelastic diffusion materials based on  Green and Naghdi
theory. In the present paper we use the results of \cite{Aouadi} to
derive a bending theory  of thermoelastic diffusion thin plates of
Mindlin type  in the context of the
 Green and Naghdi theory (models II and III). The theory allows for the effect of transverse shear
deformation as in the Mindlin-Timoshenko model of plates (see
\cite{Lagnese}), but we do not introduce any shear correction
factor.

 A linear theory of thermoelastic plates with voids was
investigated by \cite{Birsan}  under the classical theory of
thermoelasticity based on Fourier's law. 
\cite{Quintanilla} presented the basic equations which characterize
the bending of thin microstretch elastic plates under Fourier's law.
\cite{Leseduarte}
derived a bending theory  of thermoelastic  plates  in the context
of  Green and Naghdi's theory (model III). In \citep{Ghiba1, Ghiba2},
Ghiba studies the temporal and the spatial behaviour of the solution
of the bending plates of Mindlin type thermoelastic  with voids.

The organization of this paper is as follows. In Section 2 we
describe the theory established by \cite{Aouadi} to
obtain, in section 3,  the bending theory of thermoelastic diffusion
plates based on Green-Naghdi theory of type II and III, which admits
the possibility of "second sound".
 With the help of the semigroup
theory of linear operators  we investigate the well-posedness and
the asymptotic behavior  of solutions to the proposed model in
sections 4 and 5. In section 6 we introduce a weighted surface
measure associated with the dynamic process at issue  and then we
establish a second-order differential inequality whose integration
gives a good information upon the spatial behavior. Finally,  in
Section 7, we study the problem of localization in time of the
solutions. For this end, we use the uniqueness property of the
backward in time problem.

 It is
worth noting that we focus on the analysis of the qualitative
properties of solutions of type III problem. However, some
particular aspects of the type II problem are also pointed out.

\section{ Basic Equations and Preliminaries}

We refer the motion of the continuum to a fixed system of
rectangular Cartesian axes $Ox_i\ (i = 1,2, 3)$. We shall employ the
usual summation and differentiation  conventions. Latin subscripts
are understood to range over the integers $(1, 2, 3)$ summation over
repeated subscripts is implied and subscripts preceded by a comma
denote partial differentiation with respect to the corresponding
material Cartesian coordinate. We consider a body that at time $t_0$
occupies a bounded  regular region $V$ of the Euclidean
three-dimensional space and is bounded by the surface $\p V$. We
deal with functions of position and time which have as their domain
of definition the Cartesian product $\bar V\times [0,\infty) $,
where $\bar V$ is the closure of $V$. Letters in boldface stand for
tensors of order $p\geq 1$ and, if $\bf{v}$ has order $p$, we write
$v_{i_1,i_2,\cdots,i_p}$ in the Cartesian coordinate frame. In what
follows, we use a superposed dot to denote material time
differentiation.

On the basis of the  theory established by 
\cite{Aouadi}, the behavior of thermoelastic diffusion bodies  under
Green-Naghdi theory  is governed by the following equations:

\begin{itemize}
  \item [ (i)] the motion equation
 \be\label{motion} t_{ji,j}+ f^*_i=\rho
\ddot u_i,\ee  \item [ (ii)] the energy equation \be\label{energy} \rho  \dot
S=-\Phi_{i,i}+ s^*,\ee  \item [ (iii)] the equation of conservation of mass diffusion  \be\label{mass} \dot
C=-\eta_{i,i}+ c^*.\ee
\end{itemize}
If the material is isotropic, then the constitutive equations
become
\begin{align}\begin{split}\label{const}
 t_{ij} &=\lambda e_{kk}\delta_{ij}+2\mu e_{ij}-d_{1}T\delta_{ij}-d_{2}\psi\delta_{ij},\\
 \rho  S&=\displaystyle{d_{1} e_{kk}+cT+\kappa  \psi},\\
  C&=d_{2} e_{kk}+\kappa T+r \psi,\\
   \Phi_i&=-(k_1 \nu+\hbar_1 \gamma+k_2 T+\hbar_2 \psi)_{,i},\\
   \eta_i&= -(h_1 \gamma+\hbar_1 \nu+\hbar_2T+h_2 \psi)_{,i},   \\
   q_i&=T_0\Phi_i=-T_0(k_1 \nu+\hbar_1 \gamma+k_2 T+\hbar_2 \psi)_{,i}.\\
\end{split}\end{align}

In these equations we have used  the following notations:
$t_{ij}$ is the stress tensor, $f^*_i$ is the body force per
unit of initial volume, $\rho$  is   the reference mass density,
$T$ is the absolute temperature, $\textbf{u}=(u_i)$  is
the displacement vector,  $\Phi_i$ is the entropy flux vector,  $q_i$ is the heat flux vector, $\eta_i$ denotes the flow of the diffusing mass vector,
 $s^*$  is the external  rate of supply of entropy per unit mass,
 $e_{ij}=  \half(u_{i,j}+u_{j,i})$  is  the deformation tensor,  $S$ is the entropy per unit
mass,  $\psi$ is the chemical potential per unit mass,  $C$ is the
concentration of diffusive material in the elastic body, $\nu$ is
the thermal displacement whose derivative coincides with the
absolute temperature, $\psi$ is the chemical potential displacement
whose derivative coincides with the chemical potential,
i.e.,\[\nu=\int_{t_0}^t T dt,\ \ \ \ \gamma=\int_{t_0}^t \psi dt.\]
These scalars, on the macroscopic scale, are regarded, respectively,
as representing some "mean" thermal and chemical potential
displacements magnitudes on the molecular scale. The constants
$\lambda$ and $\mu$ are  elastic coefficients and $d_1$ and $d_2$
are the coefficients of thermal and mass diffusion expansions,
respectively. The constants $\kappa$ and $r$ are measures of
thermodiffusion and diffusive effects, respectively. $k_1$ and $h_1$
are coefficients of thermal and diffusion conductivity,
respectively, while $\hbar_1$ is a measure of thermodiffusion
gradient displacement.  $k_2$, $h_2$ and $\hbar_2$   are
coefficients characterizing the type III model.

Remark that the evolutive equations for the thermoelastic diffusion theory of type II (without energy dissipation) can
be deduced from the above equations
by taking $k_2=h_2=\hbar_2=0$.

All the above coefficients are constitutive constants and satisfy
the following conditions :

\begin{itemize}
  \item [ (i)]  To stabilize the thermoelastic diffusion system, we need
  that \citep{Aouadi}
 \be\label{cond1}
   \delta=cr-\kappa^2>0,
 \ee
 which implies that \[ c\theta^2+2\kappa \theta P+rP^2>0.\]

  \item [ (ii)] For type III model, the coefficients $k_2$, $h_2$ and $\hbar_2$  should satisfy
 \be\label{cond2}
   k_2h_2-\hbar_2^2>0,
 \ee
  which implies that \[k_2\theta^2+2\hbar_2 \theta P+h_2P^2>0\] to ensure the
non-negativeness of the  internal rate of production of entropy (see \cite{Aouadi}).

\item [ (iii)] For both type II and III models,  the necessary
and sufficient conditions that the internal energy density be
positive are \be\label{cond3} \mu>0,\ \ \  \lambda+\mu>0,\ \ \   k_1h_1-\hbar_1^2>0. \ee
\end{itemize}

The components of surface traction $t_i$,  the internal flux of
entropy per unit mass $\Phi$ and  the diffusion flux $\eta$  at
regular points of $\p V$, are given by \be t_i=t_{ji}n_j,\ \ \ \
\Phi=\Phi_i n_i,\ \ \ \ \eta=\eta_i n_i,\ee respectively. We denote
by $n_j$ the outward unit normal of $\p V$. We assume that $(f_i^*,\
s^*,\ c^*)$ are continuous on $V \times \mathcal{T}$, where
$\mathcal{T}$ is a temporal interval. To the system of field
equation we must add boundary conditions and initial conditions.

\section{Bending  thermoelastic diffusion plates}

In what follows we assume that the region $V$ is the interior of a
right cylinder of length $2 h$ with open cross-section $\Sigma$ and
smooth lateral boundary $\Pi$. Let $\Gamma$ be the boundary of
$\Sigma$. The rectangular Cartesian coordinate frame is supposed to
be chosen in such a way that the plane $x_1Ox_2$ is the middle
plane. Thus, we have
\[V=\{{\bf x}:(x_1,x_2)\in \Sigma,\ -h<x_3<h\},\ \ \ \Pi=\{{\bf x}:(x_1,x_2)\in \Gamma,\ -h<x_3<h\}.\]
An
admissible process $p=\{u_i,\nu,T,\gamma,\psi\}$ is a state of bending
on $V\times\mathcal{T}$ provided
\begin{align}
  \begin{split}\label{sym}
    u_\alpha (x_1,x_2,x_3,t)  &=- u_\alpha(x_1,x_2,-x_3,t),\qquad
      u_3 (x_1,x_2,x_3,t)  = u_3 (x_1,x_2,-x_3,t),\\
    \nu (x_1,x_2,x_3,t) &= -\nu (x_1,x_2,-x_3,t),\qquad
      \gamma(x_1,x_2,x_3,t)  =- \gamma (x_1,x_2,-x_3,t),\quad(x_1,x_2,x_3,t) \in V \times \mathcal{T}.
  \end{split}
\end{align}
Here and in what follows the Greek subscripts are confined to the
range $1,\ 2$. In view of Eqs. \eqref{motion}, \eqref{const} we find that
\begin{align}
\begin{split}\label{sym2}
t_{\alpha \beta}(x_1,x_2,x_3,t)  &=- t_{\alpha \beta}(x_1,x_2,-x_3,t),\qquad
t_{33}(x_1,x_2,x_3,t) =- t_{33}(x_1,x_2,-x_3,t),\\
t_{\alpha 3}(x_1,x_2,x_3,t)  &= t_{\alpha 3} (x_1,x_2,-x_3,t),\qquad
t_{3\alpha }(x_1,x_2,x_3,t)  = t_{3\alpha } (x_1,x_2,-x_3,t),\\
S(x_1,x_2,x_3,t) &=- S (x_1,x_2,-x_3,t),\qquad
\Phi_{\alpha}(x_1,x_2,x_3,t)  =- \Phi_{\alpha } (x_1,x_2,-x_3,t),\\
\Phi_{3}(x_1,x_2,x_3,t)  &=\Phi_{3 } (x_1,x_2,-x_3,t),\qquad
q_{\alpha}(x_1,x_2,x_3,t)  =- q_{\alpha } (x_1,x_2,-x_3,t),\\
q_{3}(x_1,x_2,x_3,t)  &=q_{3 } (x_1,x_2,-x_3,t),\qquad
C(x_1,x_2,x_3,t)=- C (x_1,x_2,-x_3,t),\\
\eta_{\alpha }(x_1,x_2,x_3,t)  &=- \eta_{\alpha } (x_1,x_2,-x_3,t),\qquad
\eta_{3 }(x_1,x_2,x_3,t)=\eta_{3 } (x_1,x_2,-x_3,t).
\end{split}
\end{align}
We say that the
system of body loads $(f_i^*,\ s^*,\ c^*)$ is compatible with a
state of bending if
\begin{align}
\begin{split}\label{rest}
f^*_{\alpha }(x_1,x_2,x_3,t)  &=- f^*_{\alpha }(x_1,x_2,-x_3,t),\qquad
f^*_{3 }(x_1,x_2,x_3,t) = f^*_{3 }(x_1,x_2,-x_3,t),\\
s^*(x_1,x_2,x_3,t) &=- s^*(x_1,x_2,-x_3,t),\qquad
c^*(x_1,x_2,x_3,t) =- c^*(x_1,x_2,-x_3,t).
\end{split}
\end{align}
In what follows we assume
that the body loads satisfy the restrictions \eqref{rest}. From Eqs.
\eqref{sym2}, \eqref{rest} we get
\begin{align}
\begin{split}\label{nul}
\int_{-h}^{h} t_{\alpha \beta} dx_3&=0,\quad
\int_{-h}^{h}t_{33}dx_3 =0,\quad
\int_{-h}^{h} Sdx_3 =0,\\
\int_{-h}^{h} \Phi_{\alpha } dx_3&=0,\quad
\int_{-h}^{h} q_{\alpha } dx_3=0,\quad
\int_{-h}^{h}C dx_3=0,\\
\int_{-h}^{h} \eta_{\alpha }dx_3 &=0,\quad
\int_{-h}^{h} f^*_{\alpha } dx_3=0,\quad
\int_{-h}^{h} s^* dx_3 =0,\quad
\int_{-h}^{h} c^* dx_3=0.
\end{split}
\end{align}
We derive a theory of thin plates of uniform thickness assuming that
the fields $u_i$, $\nu$ and $\gamma$ do not vary violently with
respect to $x_3$. We denote \be\label{defN} N_\al=\frac{1}{2h}
\int_{-h}^{h} t_{\alpha 3} dx_3. \ee We assume that the functions
$t_i,$  $q$ and $\eta$ are prescribed on the surfaces $x_3=\pm h$.
We integrate \eqref{motion} with respect to $x_3$ between the limits
$-h$ and $h$. According to Eqs. \eqref{sym} and \eqref{defN}, we
obtain \be\label{motion1} N_{\al,\al}+f=\rho\ddot w,\quad\text{on }
\Sigma\times\mathcal{T}, \ee where \be w=\frac{1}{2h} \int_{-h}^{h}
u_{3} dx_3,\quad f=\frac{1}{h}t_{33}(x_1,x_2,h,t)+ \frac{1}{2h}
\int_{-h}^{h} f^*_{3} dx_3. \ee If we multiply \eqref{motion}  by
$x_3$  and integrate from $x_3=-h$ to $x_3=h$, then we obtain
\be\label{motion2} M_{\beta\alpha,\beta}-2hN_\al+f_\al=\rho I \ddot
v_\al, \quad\text{on }\Sigma\times \mathcal{T}, \ee where
\begin{align}\begin{split}\label{defM} M_{\al\beta}&= \int_{-h}^{h} x_3t_{\al\beta} dx_3,
\quad
 I v_{\al}=
\int_{-h}^{h} x_3u_{\al} dx_3,\quad I= \frac{2}{3}h^3,\cr
f_\al&=2ht_{3\alpha}(x_1,x_2,h,t)+ \int_{-h}^{h} x_{3}f^*_\al
dx_3.
\end{split}\end{align}
Now multiply  \eqref{energy} and \eqref{mass} by $x_3$ and integrate
from $x_3=-h$ to $x_3=h$. In view of \eqref{sym}-\eqref{rest} we
obtain, respectively \be\label{temp} \rho \dot
\sigma=-\Psi_{\alpha,\alpha}+2hR+ W,\ee\be\label{diff} \dot
\chi=-\Omega_{\alpha,\alpha}+2hM+V,\ee

where
\begin{align}
\begin{split}\label{diftempmass}
\sigma&=\int_{-h}^{h} x_3S dx_3,\quad \Psi_\alpha=
\int_{-h}^{h}x_3\Phi_{\alpha} dx_3,\quad
R= \frac{1}{2h}\int_{-h}^{h} \Phi_3dx_3,\\
W&=- 2h\Phi_{3}(x_1,x_2,h,t)+ \int_{-h}^{h} x_{3}s^* dx_3,\\
\chi&= \int_{-h}^{h} x_3C dx_3,\quad
\Omega_\alpha= \int_{-h}^{h}x_3\eta_{\alpha} dx_3,\quad
M= \frac{1}{2h}\int_{-h}^{h} \eta_3dx_3,\\
V&=- 2h\eta_{3}(x_1,x_2,h,t)+ \int_{-h}^{h} x_{3}c^*dx_3.
\end{split}
\end{align}
The functions $f,\ f_\al,\ W$ and $V$
are prescribed.

We restrict our attention to the state of
bending characterized by
\begin{align}\begin{split}\label{new} u_\alpha&= x_3v_\alpha (x_1,x_2,t),\qquad
 u_3= w (x_1,x_2,t),\cr
\nu &= x_3\tau(x_1,x_2,t),\qquad T =x_3\theta(x_1,x_2,t),\cr \gamma
&= x_3\wp(x_1,x_2,t),\qquad \psi = x_3P(x_1,x_2,t),\quad\text{on }
\Sigma\times \cal{T}.\end{split}\end{align}

In view of \eqref{new} we have
\begin{align}\begin{split}\label{new2}
e_{\al\beta}&=x_3\varepsilon_{\al\beta},\ \ \
2e_{\al3}=\gamma_{\al},\ \ \ e_{33}=0,\cr
\varepsilon_{\al\beta}&=\half(v_{\al,\beta}+v_{\beta,\alpha}),\ \ \
\gamma_\al=v_\al+w_{,\al}=2\varepsilon_{3\alpha}=2\varepsilon_{\al3}.\end{split}\end{align}

The quantities $\gamma_\al$ represent the angles of rotation of
the cross sections $x_\al = \text{const}$ about the middle surface.

It follows from \eqref{const}, \eqref{new} and \eqref{new2} that
\begin{align}
\begin{split}\label{constv}
t_{\alpha\beta}&=x_3(\lambda\varepsilon_{rr}\delta_{\alpha\beta}+2\mu\varepsilon_{\alpha\beta}-
  d_1\theta\delta_{\alpha\beta}-d_2 P\delta_{\alpha\beta}),\\
t_{\alpha 3}&=\mu \gamma_\al,\\
t_{33}&=x_3(\lambda \varepsilon_{rr}-d_1 \theta-d_2 P),\\
\rho S&=x_3(d_1 \varepsilon_{rr}+c \theta+\kappa P),\\
C&=x_3(d_2 \varepsilon_{rr}+\kappa \theta+r P),\\
\Phi_\alpha&=-x_3(k_1 \tau+\hbar_1 \wp+k_2 \theta+\hbar_2 P)_{,\alpha},\\
\Phi_{3}&=-(k_1 \tau+\hbar_1 \wp+k_2 \theta+\hbar_2 P),\\
\eta_\alpha&= -x_3(h_1 \wp+\hbar_1 \tau+\hbar_2\theta+h_2 P)_{,\alpha},\\
\eta_3&= -(h_1 \wp+\hbar_1 \tau+\hbar_2\theta+h_2 P),\\
q_\alpha&=-x_3T_0(k_1 \tau+\hbar_1 \wp+k_2 \theta+\hbar_2 P)_{,\alpha},\\
q_{3}&=-T_0(k_1 \tau+\hbar_1 \wp+k_2 \theta+\hbar_2 P).
\end{split}
\end{align}

It follows from Eqs. \eqref{defN}, \eqref{defM}, \eqref{diftempmass} and \eqref{constv} that
\begin{align}\begin{split}\label{constvv} M_{\alpha\beta}&=I(\lambda
\varepsilon_{rr}\delta_{\alpha\beta}+2\mu
\varepsilon_{\alpha\beta}-d_1\theta\delta_{\alpha\beta}-d_2 P\delta_{\alpha\beta}),\cr
 N_{\alpha }&=\mu \gamma_\al,\cr
\rho\sigma&=I(d_1 \varepsilon_{rr}+c \theta+\kappa P),\cr
\Psi_{\alpha}&=-I(k_1 \tau+\hbar_1 \wp+k_2 \theta+\hbar_2
P)_{,\alpha},\cr R&=-(k_1 \tau+\hbar_1 \wp+k_2 \theta+\hbar_2 P),\cr
\chi&=I(d_2 \varepsilon_{rr}+\kappa \theta+r P),\cr
\Omega_{\alpha}&=-I(h_1 \wp+\hbar_1 \tau+\hbar_2\theta+h_2
P)_{,\alpha},\cr M&=-(h_1 \wp+\hbar_1 \tau+\hbar_2\theta+h_2
P).\end{split}\end{align}

The equations of thermoelastic diffusion bending plates  consist of
the equations of motion \eqref{motion1} and \eqref{motion2}, the
entropy equation \eqref{temp}, the mass diffusion equation
\eqref{diff}, the constitutive equations \eqref{constvv} and the
geometrical equations \eqref{new2}. These equations can be expressed
in terms of the functions $v_\alpha,w,\tau,\theta,\wp$ and $P$.
Thus, we obtain the evolutive equations
   of type III  (with energy dissipation)
\begin{align}\begin{split}\label{equat} \rho I
\ddot v_\alpha&=I(\mu \Delta
v_\alpha+(\lambda+\mu)v_{\beta,\beta\alpha}-d_1
\theta_{,\alpha}-d_2 P_{,\alpha})-2h\mu(v_\alpha+w_{,\alpha})+f_\alpha,\cr  \rho
\ddot w
 &=\mu \Delta w+\mu v_{\alpha,\alpha}+f,\cr cI\ddot \tau+\kappa I\ddot \wp&=
 I(k_1 \Delta\tau+\hbar_1 \Delta\wp+k_2 \Delta\theta+\hbar_2 \Delta P)-Id_1 \dot v_{\alpha,\alpha}-2h(k_1 \tau+\hbar_1 \wp+k_2 \theta+\hbar_2 P)+ W,\cr
\kappa I\ddot \tau+r I\ddot \wp&= I(h_1 \Delta\wp+\hbar_1 \Delta\tau+\hbar_2\Delta\theta+h_2\Delta P)-Id_2 \dot v_{\alpha,\alpha}-2h(h_1 \wp+\hbar_1 \tau+\hbar_2\theta+h_2 P)+ V.\end{split}\end{align}

Remark that the evolutive equations  of the thermoelastic diffusion bending plates
 of type II (without energy dissipation) can be deduced from   (\ref{equat})  by taking $k_{2}=h_{2}=\hbar_{2}=0$.

Summarizing, the following initial boundary value problems
are to be solved:

$(\verb"i")$  Type II problem : Find $(v_\alpha,z_\alpha,w,y,
\tau,\theta,\wp,P)$ solution of (\ref{equat}) (with
$k_{2}=h_{2}=\hbar_{2}=0$) subject to the initial conditions
\begin{align}
\begin{split}\label{initial}
v_\alpha(x_1,x_2,0)&=v_\alpha^0(x_1,x_2),\quad\dot v_\alpha(x_1,x_2,0)=z^0_\alpha(x_1,x_2),\\
w(x_1,x_2,0)&=w^0(x_1,x_2),\quad\dot w(x_1,x_2,0)= y^0(x_1,x_2),\\
\tau(x_1,x_2,0)&=\tau^0(x_1,x_2),\quad\theta(x_1,x_2,0)= \theta^0(x_1,x_2),\\
\wp(x_1,x_2,0)&=\wp^0(x_1,x_2),\quad P(x_1,x_2,0)=P^0(x_1,x_2), \quad (x_1,x_2) \in \bar \Sigma,
\end{split}
\end{align}
where $v_\alpha^0,\ z_\alpha^0,\ w^0,\ y^0,\  \tau^0,\  \theta^0,\  \wp^0$ and $P^0$  are
prescribed functions.

We consider the boundary conditions
\be\label{boundary}
M_{\beta\alpha}n_\beta =\tilde{M}_\alpha,\quad
N_{\alpha}n_\alpha =\tilde{N},\quad
\Psi_\alpha n_\alpha =\tilde{\Psi},\quad
\Omega_\alpha n_\alpha =\tilde{\Omega},\quad\text{on }\Gamma\times \mathcal{T},
\ee
where the given
functions $\tilde{M}_\alpha,\ \tilde{N},\ \tilde{\Psi}$ and
$\tilde{\Omega}$ are piecewise regular and continuous in time.

$(\verb"ii")$ Type III problem : Find $(v_\alpha,z_\alpha,w,y, \tau,\theta,\wp,P)$  solution of (\ref{equat}) subject to the initial conditions (\ref{initial})
and the boundary conditions (\ref{boundary}).

Some qualitative properties of  the solutions of type III problem
are studied in the following.
However, some particular aspects of the type II problem are also pointed out.

\section{Well-posedness}
We shall use the results of the semigroup of linear operators theory
to  prove the existence, uniqueness and continuous dependence from
the initial values and the external loads of the solution for the
system  (\ref{equat}).  Seeking for simplicity, we will restrict
ourselves to homogeneous boundary conditions \be
\label{BC}{{v}_\alpha}=0,\ \ \ \ w=0,\ \ \ \tau=0,\ \ \ \ \wp=0
\quad  \hbox{on}\ \ \Gamma \times \mathcal{T}.\ee

In the rest of the paper we assume:
\begin{itemize}
  \item [ (i)] Relations (\ref{cond1})-(\ref{cond3}) are satisfied;
\item [ (ii)]  The positive definiteness of the internal energy density, i.e., there exists a positive constant $c_0$ such that
\end{itemize} \begin{align}\begin{split}\label{positive}&I(
 \lambda \varepsilon_{rr}\varepsilon_{\gamma\gamma}+ 2\mu
 \varepsilon_{\alpha\beta}\varepsilon_{\alpha\beta})+2h\mu\gamma_\alpha\gamma_\alpha
+I(k_1\tau_{,\alpha}\tau_{,\alpha}+h_1\wp_{,\alpha}\wp_{,\alpha}+2\hbar_1\tau_{,\alpha}\wp_{,\alpha})
+2h(k_1\tau^2+h_1\wp^2+2\hbar_1\tau\wp)\cr&\ge c_0(
\varepsilon_{\alpha\beta}\varepsilon_{\alpha\beta}+\gamma_\al
\gamma_\al+\tau_{,\alpha}\tau_{,\alpha}+\wp_{,\alpha}\wp_{,\alpha}+\tau^2+\wp^2
)\end{split}\end{align} for any $\varepsilon_{\alpha
\beta},\gamma_\al,\tau,\tau_{,\alpha },\wp,\wp_{,\alpha }$.

 We now wish to transform the boundary-initial-value
problem defined by system (\ref{equat}), the initial conditions
(\ref{initial}) and the boundary conditions (\ref{BC}) to an
abstract problem on a suitable Hilbert space. In what follows we use
the notation $\textbf{z}=\dot {\textbf{v}},\ y=\dot {w},\
\theta=\dot \tau,\ P=\dot \wp$. Let
 \[\mathscr{H}=\left\{(\textbf{v}, \textbf{z}, {w}, {y},\tau,\theta,
\wp,P);\
  {v}_\alpha, w,\tau,\wp\in  {{W}}_0^{1,2}(\Sigma);\ {z}_\alpha, y,\theta,P\in  {{L}}^{2}(\Sigma)
 \right\},\] where ${W}_0^{1,2}(\Sigma)$ and $ L^{2}(\Sigma)$ are
 the familiar Sobolev spaces.

We consider the following
operators\begin{align}\begin{split}\label{matrix}  A_\alpha
\textbf{v}  &= \rho^{-1}[\mu \Delta
v_\alpha+(\lambda+\mu)v_{\beta,\beta\alpha}]-2h\rho^{-1}I^{-1}\mu
v_\alpha,\qquad B_\alpha w=-2\rho^{-1}I^{-1}h\mu w_{,\alpha},\cr
C_\alpha \theta&=-\rho^{-1}d_1\theta_{,\alpha},\qquad D_\alpha P=
-\rho^{-1}d_2P_{,\alpha},\qquad G \textbf{v}  = \rho^{-1}\mu
v_{\alpha,\alpha},\qquad
 \mathcal{H} w=\rho^{-1}\mu\Delta w,\cr
{J} \textbf{z}  &= -\delta^{-1} (rd_1 -\kappa d_2)
z_{\alpha,\alpha},\qquad
 \mathcal{K}\tau = \delta^{-1}(rk_1 -\kappa \hbar_1)\Delta \tau-2h (I\delta)^{-1}(rk_1 -\kappa \hbar_1)\tau,\cr
 \mathcal{L}\theta &=\delta^{-1} (rk_2-\kappa \hbar_2) \Delta\theta -2h (I\delta)^{-1} (rk_2-\kappa \hbar_2)
 \theta,\qquad
 \mathcal{M}\wp = \delta^{-1}(r \hbar_1-\kappa h_1)\Delta \wp-2h (I\delta)^{-1}(r \hbar_1-\kappa h_1) \wp,\cr
 \mathcal{N} P  &=\delta^{-1}(r \hbar_2-\kappa h_2)\Delta P-2h (I\delta)^{-1}(r \hbar_2-\kappa h_2)
 P,\qquad
  L  \textbf{z} =-\delta^{-1} (cd_2 -\kappa d_1) z_{\alpha,\alpha},\cr
 \mathcal{U}\tau &= \delta^{-1}(c\hbar_1 -\kappa k_1)\Delta \tau-2h (I\delta)^{-1}(c\hbar_1 -\kappa
 k_1)\tau,\qquad
 \mathcal{V}\theta = \delta^{-1} (c\hbar_2-\kappa k_2) \Delta\theta-2h (I\delta)^{-1} (c\hbar_2-\kappa k_2) \theta,\cr
 \mathcal{W}\wp &= \delta^{-1}(c h_1-\kappa \hbar_1)\Delta \wp-2h (I\delta)^{-1}(ch_1-\kappa \hbar_1) \wp,\qquad \mathcal{X} P  =
\delta^{-1}(c h_2-\kappa \hbar_2)\Delta P-2h (I\delta)^{-1}(c
h_2-\kappa \hbar_2) P,
\end{split}\end{align} where $\delta$ is given by \eqref{cond1}, $\textbf{A}  =(A_i),\ \textbf{B}  =(B_i),\
\textbf{C} =(C_i)$ and
 $\textbf{D}  =(D_i)$.
Now consider the matrix operator $\mathscr{A} $ on $\mathscr{H}$
defined by \[\mathscr{A}= \left( {\begin{array}{*{20}c}
  0& {\bf Id} &    0&  0&  0&  0&0&  0&   \\
      {\bf A}&  0&{\bf B}&  0 & 0&{\bf C}  &0& {\bf D} \\
0&0&0&{ Id} &    0&   0&0&  0&  \\
{ G}&  0& \mathcal{H}&   0&0  & 0&0&  0&\\
0&0&0&0 &      0& { Id}&0&  0& \\
 0&  {J}& 0&   0&\mathcal{K}  & \mathcal{L}&\mathcal{M}&  \mathcal{N}&\\
 0&0&0&0 &      0& 0&0&  Id& \\
0&  {L}& 0&   0&\mathcal{U}  & \mathcal{V}&\mathcal{W}&  \mathcal{X}&\\
\end{array}}\right)\]

with the domain \begin{align}\begin{split}\mathscr{D} =\mathscr{D}
(\mathscr{A}
)&=({\textbf{W}}_0^{1,2}\cap{\textbf{W}}^{2,2})\times{\textbf{W}}_0^{1,2}\times({{W}}_0^{1,2}\cap{{W}}^{2,2})\times{{W}}_0^{1,2}
\times({{W}}_0^{1,2}\cap{{W}}^{2,2})\times({{W}}_0^{1,2}\cap{{W}}^{2,2})\cr&\times({{W}}_0^{1,2}\cap{{W}}^{2,2})\times({{W}}_0^{1,2}\cap{{W}}^{2,2}),\end{split}\end{align}
where $\textbf{Id}$ and $Id$ are  the identity operators in the
respective spaces. We note that the domain $\mathscr{D}$ is dense in
$\mathscr{H}$.

 We
introduce the inner product in $\mathscr{H}$ defined by
\begin{align}\begin{split}\label{ps} \langle(\textbf{v}, \textbf{z},
{w}, {y},\tau,\theta,\wp,P ),(\textbf{v}^*, \textbf{z}^*, {w}^*,
{y}^*,\tau^*,\theta^*, \wp^*,P^*)\rangle &=\half\int_\Sigma\Big(\rho
I {z_\alpha}{z_\alpha}^* +2h\rho yy^*+c I \theta\theta^*+\kappa I
(\theta^* P+\theta P^*)\cr& +rI PP^*+\mathscr{W}[(\textbf{v}, {w},
\tau, \wp),(\textbf{v}^*, {w}^*, \tau^*, \wp^*)] \Big)da,
\end{split}\end{align}  where
\[\mathscr{W}[(\textbf{v}, {w}, \tau, \wp),(\textbf{v}^*, {w}^*,
\tau^*, \wp^*)]= I\Big(
\lambda\varepsilon_{\alpha\alpha}\varepsilon_{\beta\beta}^*
+2\mu\varepsilon_{\alpha\beta}\varepsilon^*_{\alpha\beta}+k_1\tau_{,\alpha}\tau_{,\alpha}^*+h_1\wp_{,\alpha}\wp_{,\alpha}^*+
\hbar_1(\tau_{,\alpha}^*\wp_{,\alpha}+\tau_{,\alpha}\wp_{,\alpha}^*)\Big)\]
\[~~~~~~~~~~~~~+2h\Big(\mu
\gamma_\alpha\gamma_\alpha^*+k_1\tau\tau^*+h_1\wp\wp^*+
\hbar_1(\tau^*\wp+\tau\wp^*)\Big).\]

In the above relations we have used the notation
\[ \varepsilon_{\alpha\beta}^*=\half(v^*_{\beta,\alpha}+v^*_{\alpha,\beta}),\ \ \ \ \gamma_\al^*=v^*_{\alpha}+w^*_{,\alpha}
=2\varepsilon^*_{\alpha 3}=2\varepsilon^*_{3\alpha}.  \]

If we recall the assumptions (\ref{positive}) and the first Korn inequality,
 we conclude that
the norm induced in $\mathscr{H}$ through the product (\ref{ps}) is
equivalent to the usual one in $\mathscr{H}$.

The boundary initial value problem (\ref{equat}),  (\ref{initial}),
(\ref{BC}) can be transformed into the following equation in the
Hilbert space $\mathscr{H}$, \be\label{cauchy}
\frac{dU(t)}{dt}={\mathscr{A}} U(t)+{\mathcal{F}}(t),\ \ \
U(0)=U_0,\ee where
\[U=(\textbf{v}, \textbf{z}, {w}, {y},\tau,\theta,\wp,P),\ \ \ \ U_0=(v_\alpha^0, z_\alpha^0, {w}^0,
{y}^0,\tau^0,\theta^0, \wp^0,P^0),\]
\[{\mathcal{F}}=(0,(\rho I)^{-1}f_\al,0,
\rho^{-1} f,0,(I\delta)^{-1}(rW -\kappa V),0,(I\delta)^{-1}(cV
-\kappa W) ).\]

Now, we use the theory of semigroups of linear operators to obtain
the existence of solutions for the Eq. (\ref{cauchy}).

\begin{lemma}  The operator $\mathscr{A}$ satisfies the inequality
$<{\mathscr{A}}U,U>< 0$ for type III model and the equality
$<{\mathscr{A}}U,U>= 0$ for type II model, for every $U \in
\mathscr{D}(\mathscr{A})$.
\end{lemma}

{\bf Proof.} Let $U=(\textbf{v}, \textbf{z}, {w}, {y},\tau,\theta,
\wp,P)\in\mathscr{D}(\mathscr{A})$. Using the divergence theorem and
the boundary conditions, we have
\begin{align}\begin{split}\label{dissp}<{\mathscr{A}}U,U>&=-
\int_\Sigma\Big(z_{\alpha,\beta}M_{\beta\alpha}+2hz_{\alpha}N_{\alpha}+2hy_{,\alpha}N_{\alpha}+I(k_2\theta_{,\alpha}\theta_{,\alpha}
+2\hbar_2\theta_{,\alpha}
P_{,\alpha}+h_2P_{,\alpha}P_{,\alpha})\cr&+2h(k_2\theta^2
+2\hbar_2\theta P+h_2P^2) +Id_1z_{r,r}\theta+Id_2z_{r,r}P-
{\mathscr{W}}[(\textbf{v}, {w}, \tau, \wp),(\textbf{v}, {w}, \tau,
\wp)]\Big)da\cr&=-\int_\Sigma
\Big(I(k_2\theta_{,\alpha}\theta_{,\alpha} +2\hbar_2\theta_{,\alpha}
P_{,\alpha}+h_2P_{,\alpha}P_{,\alpha})+2h(k_2\theta^2
+2\hbar_2\theta P+h_2P^2)\Big)da.\end{split}\end{align}

In the context  of type III model and from the assumption
(\ref{cond2}) we have $<{\mathscr{A}}{U},{U}><0$, which means
dissipation of the energy, i.e.,
$$\frac{dE_0(t)}{dt}=-\int_\Sigma
\Big(I(k_2\theta_{,\alpha}\theta_{,\alpha} +2\hbar_2\theta_{,\alpha}
P_{,\alpha}+h_2P_{,\alpha}P_{,\alpha})+2h(k_2\theta^2
+2\hbar_2\theta P+h_2P^2)\Big)da<0,$$ where $$
E_0(t)=\half\int_\Omega\Big(\rho I {z_\alpha}{z_\alpha} +2h\rho yy+c
I \theta^2+2\kappa I \theta P+rI P^2 +2\mathscr{V}  \Big)da,$$ and
\be \label{V}2\mathscr{V}= I\Big(
\lambda\varepsilon_{\alpha\alpha}\varepsilon_{\beta\beta}
+2\mu\varepsilon_{\alpha\beta}\varepsilon_{\alpha\beta}+k_1\tau_{,\alpha}\tau_{,\alpha}+2\hbar_1\tau_{,\alpha}\wp_{,\alpha}+h_1\wp_{,\alpha}\wp_{,\alpha}
\Big)+2h\Big(\mu \gamma_\alpha\gamma_\alpha+k_1\tau^2+ 2\hbar_1
\tau\wp+h_1\wp^2\Big).\ee
 In the context  of
type II model ($ k_2=h_2=\hbar_2 = 0$), we have
$<{\mathscr{A}}{U},{U}>=0$, which means conservation of the energy
i.e., $\frac{dE_0(t)}{dt}=0$.\hfill$\Box$

It is worth remarking that this quantity is also conserved even if
we do not impose conditions $(\textrm{i})-(\textrm{ii})$.

\begin{lemma}  The operator $\mathscr{A}$ has the property that \[\hbox{Range}(\mathcal{I }-\mathscr{A})=\mathscr{H},\]
where $\mathcal{I}$ is the identity operator in $\mathscr{H}$.
\end{lemma}

{\bf Proof:} Let $U^*=(\textbf{v}^*, \textbf{z}^*, {w}^*,
{y}^*,\tau^*,\theta^*, \wp^*,P^*)\in \mathscr{H}$. We must prove
that
\[U-{\mathscr{A}}U =U^*\] has a solution
$U=(\textbf{v}, \textbf{z}, {w}, {y},\tau,\theta,\wp,P )$
in $\mathscr{D}$. This equation leads to the system
 \begin{align}\begin{split}\label{surj1}
 \textbf{v}^*&=\textbf{v}-\textbf{z},\qquad w^*=w-y,\cr \tau^*&=\tau-\theta,\qquad \wp^*=\wp-P,\cr
 \textbf{z}^*&=\textbf{z}-(\textbf{Av}+\textbf{B}w+\textbf{C}\theta+\textbf{D}P),\cr
y^*&=y-(G\textbf{v}+\mathcal{H}w),\cr \theta^*
&=\theta-(J\textbf{z}+\mathcal{K}\tau+\mathcal{L}\theta+\mathcal{M}
\wp+\mathcal{N} P),\cr
P^*&=P-(L\textbf{z}+\mathcal{U}\tau+\mathcal{V}\theta+\mathcal{W}
\wp+\mathcal{X} P).\end{split}\end{align}Substituting the four first
equations in the others, we obtain
 \begin{align}\begin{split}\label{surj2}
 \textbf{z}^*+\textbf{v}^*-\textbf{C}\tau^*-\textbf{D}\wp^*&=\textbf{v}-(\textbf{Av}+\textbf{B}w+\textbf{C}\tau+\textbf{D}\wp),\cr
y^*+z^*&=w-(G\textbf{v}+\mathcal{H}w),\cr
\theta^*+\tau^*-\mathcal{L}\tau^*-\mathcal{N}\wp^*&=\tau-\Big(J\textbf{v}+(\mathcal{K}+\mathcal{L})\tau+(\mathcal{M}+\mathcal{N})\wp\Big),\cr
P^*+\wp^*-\mathcal{V}\tau^*-\mathcal{X}\wp^*&=\wp-\Big(L\textbf{v}+(\mathcal{U}+\mathcal{V})\tau+(\mathcal{W}+\mathcal{X})\wp\Big).
\end{split}\end{align}

To solve this system, we introduce the following bilinear form on
${\textbf{W}}_0^{1,2}$, \be
\mathscr{B}[(\textbf{v},w,\tau,\wp),(\hat{\textbf{v}},\hat{w},\hat{\tau},\hat{\wp})]=<({\textbf{v}}^{(1)},
{w}^{(1)},{\tau}^{(1)},{\wp}^{(1)}),(\rho  I \hat{\textbf{v}},2h
\rho \hat{w},c I \hat{\tau},r I \hat{\wp})>\ee where
 \[
\begin{array}{lll}
{\textbf{v}}^{(1)}&=&\textbf{v}-(\textbf{Av}+\textbf{B}w+\textbf{C}\tau+\textbf{D}\wp),\cr
{w}^{(1)}&=&w-(G\textbf{v}+\mathcal{H}w),\cr
{\tau}^{(1)}&=&\tau-\Big(J\textbf{v}+(\mathcal{K}+\mathcal{L})\tau+(\mathcal{M}+\mathcal{N})\wp\Big),\cr
{\wp}^{(1)}&=&\wp-\Big(L\textbf{v}+(\mathcal{U}+\mathcal{V})\tau+(\mathcal{W}+\mathcal{X})\wp\Big).\\
\end{array}\]
A direct calculation shows that $\mathscr{B}$ is bounded in each
variable. Using the divergence theorem, we have
\begin{align}\begin{split}\mathscr{B}[(\textbf{v},w,\tau,\wp),(\textbf{v},w,\tau,\wp)]
&=\int_\Sigma\Big(\rho I v_\alpha v_\alpha+2h\rho w^2+cI \theta^2+2 \kappa I\theta P+r I
P^2+\mathscr{W}[(\textbf{v},w,\tau,\wp),(\textbf{v},w,\tau,\wp)]\Big)da\\
&+\int_\Sigma \Big(I(k_2\theta_{,\alpha}\theta_{,\alpha}
+2\hbar_2\theta_{,\alpha}
P_{,\alpha}+h_2P_{,\alpha}P_{,\alpha})+2h(k_2\theta^2
+2\hbar_2\theta P+h_2P^2)\Big)da.\nonumber\end{split}\end{align} In
view of our assumptions on the constitutive coefficients, we see
that $\mathscr{B}$ is coercive on $ {\textbf{W}}^{-1,2}$. On the
other hand, it is easy to see that the vector
\[( \textbf{z}^*+\textbf{v}^*-\textbf{C}\tau^*-\textbf{D}\wp^*,y^*+z^*,\theta^*+\tau^*-\mathcal{L}\tau^*-\mathcal{N}\wp^*,P^*+\wp^*-\mathcal{V}\tau^*-\mathcal{X}\wp^*)
\] lies in $ {\textbf{W}}^{-1,2}$. Hence the
Lax-Milgram theorem  implies the existence of
$(\textbf{v},w,\tau,\wp)\in {\textbf{W}}_0^{1,2}$ which solves the
system (\ref{surj2}). Thus, the system  (\ref{surj1}) has also a
solution.\hfill$\Box$

The previous lemmas lead to next theorem.
\begin{theorem}
The operator $\mathscr{A} $ generates a semigroup of contraction in
$\mathscr{H}$.
\end{theorem}
{\bf Proof.}
The proof follows from Lumer-Phillips corollary to the
Hille-Yosida theorem \citep{Pazy}.\hfill$\Box$

 It is worth remarking that this
theorem implies that the dynamical system generated by the equations
of bending thermoelastic diffusion plates
 of type III (or type II) is stable in the sense of Lyapunov.
\begin{theorem}
Assume that $f_\alpha,\  f ,\  W ,\ V \in C^1([0,\infty),L^2)$ and
${U}_0$ is in the domain of the operator $\mathscr{A} $. Then, there
exists a unique solution ${U}(t)\in C^1([0,\infty),\mathscr{H})\cap
C^0([0,\infty),\mathscr{D})$ to the problem
(\ref{cauchy}).\end{theorem} Since the solutions are defined by
means of a semigroup of contraction, we have the estimate
\[\| {U}(t)\|  \le \| {U}_0 \|_\mathscr{H}+\int_0^t\Big(
\| {f_\alpha}(s)\|_{\bf{L}^2} +\| f(s)\|_{{L}^2}+\|
W(s)\|_{{L}^2}+\| V(s)\|_{{L}^2}\Big)ds\] which proves the
continuous dependence of the solutions upon initial data and body
loads. Thus, under assumptions $(\textrm{i})-(\textrm{ii})$ the
problem of bending thermoelastic diffusion plates
 of type III (or type II)
is well posed.

\section{Asymptotic behavior of solutions}
In this section we study the asymptotic behavior of solutions, whose
existence has been  proved previously, in the homogenous case
(${f_\alpha}$=0,\ $f$=0,\ $W$=0,\ $V$=0). In particular we are
interested in the relation between dissipation effects and time
decay of solutions. Therefore, we will continue to assume that the
assumptions $(\textrm{i})-(\textrm{ii})$ considered in the previous
section hold. However it is  worth noting that the results for this
section only hold for type III theory.

To this end, we recall that for a  semigroup of contraction,   the
precompact orbits tend to the $\omega-$limit  sets if its generator
$\mathscr{A}$ has only the fixed point ${\bf 0}$ and the structure
of the $\omega-$limit sets is determined by the eigenvectors of
eigenvalue $i \lambda$ (where $\lambda$ is a real number) in the
closed subspace
\[{\mathscr{L}}=\ll\{{U}\in \mathscr{H};\
<\mathscr{A}{U},{U}>=0\}\gg\,,
\]
where $\ll \mathcal{C} \gg$ denotes the closed vectorial subspace
generated by the set $\mathcal{C}$.

From the assumptions $\textrm{(i)-(ii)}$ it is easy to check that
$\mathscr{A}^{-1}({\textbf{0}})={\textbf{0}}$, while the
precompactness of the orbits starting in $\mathscr{D}$ is a
consequence of the following Lemma \citep{Pazy}.
\begin{lemma}  The operator $(\mathcal{I}- \mathscr{A})^{-1}$ is compact.
\end{lemma}
{\bf Proof.} Let $(\hat{\textbf{v}}_n, \hat{\textbf{z}}_n,\hat{w}_n,
\hat{y}_n,
 \hat{\tau}_n,\hat{\theta}_n, \hat{\wp}_n,\hat{P}_n)$ be a bounded sequence in
 $\mathscr{H}$ and let $({\textbf{v}}_n, {\textbf{z}}_n,{w}_n,
{y}_n, {\tau}_n, {\theta}_n,{\wp}_n, {P}_n)$ be the sequence of the
respective
 solutions of the system (\ref{surj1}). We have
 \be
\Lambda[\Gamma_n,\Gamma_n]=<\Sigma_n,\Upsilon_n>,\ee where
\[\Gamma_n=(\textbf{v}_n,w_n,\tau_n,\wp_n),\ \ \
\Sigma_n=(\textbf{v}_n^{(1)},w_n^{(1)},\tau_n^{(1)},\wp_n^{(1)}),\ \
\ \Upsilon_n=(\rho I \hat{\textbf{v}}_n, 2h \rho\hat{w}_n ,c I
\hat{\tau}_n,r I \hat{\wp}_n)\] and

 \be
\begin{array}{lll}
\textbf{v}_n^{(1)}&=&\textbf{v}_n-(\textbf{Av}_n+\textbf{B}w_n+\textbf{C}\tau_n+\textbf{D}\wp_n),\\
{w}_n^{(1)}&=&w_n-(G\textbf{v}_n+\mathcal{H}w_n),\\
\tau_n^{(1)}&=&\tau_n-\Big(J\textbf{v}_n+(\mathcal{K}+\mathcal{L})\tau_n+(\mathcal{M}+\mathcal{N})\wp_n\Big),\\
\wp_n^{(1)} &=&\wp_n-\Big(L\textbf{v}_n+(\mathcal{U}+\mathcal{V})\tau_n+(\mathcal{W}+\mathcal{X})\wp_n\Big).\\
\end{array}\ee
In view of the definition of
$(\textbf{v}_n^{(1)},w_n^{(1)},\tau_n^{(1)},\wp_n^{(1)})$, it
follows that it is  a bounded sequence in ${\bf{L}}^{2}(\Sigma)$ and
then the sequence $({\textbf{v}}_n, {w}_n,{\tau}_n,
 {\wp}_n)$ is a bounded sequence in $
{\textbf{W}}_0^{1,2}(\Sigma)$. The theorem of Rellich-Kondrasov
\citep{[25]} implies that there is  exists a subsequence converging
in ${\bf{L}}^{2}(\Sigma)$. In a similar way
\[ \textbf{z}_{n_j}=\textbf{v}_{n_j}-\hat{\textbf{v}}_{n_j},\ \ \
y_{n_j}={w}_{n_j}-\hat{w}_{n_j},\ \ \
\theta_{n_j}={\tau}_{n_j}-\hat{\tau}_{n_j},\ \ \
P_{n_j}={\wp}_{n_j}-\hat{\wp}_{n_j}\] has a sub-sequence converging
in ${\bf{L}}^{2}(\Sigma)$. Thus we conclude the existence of a
sub-sequence
\[({\textbf{v}}_{n_{j_k}},{\textbf{z}}_{n_{j_k}},
{w}_{n_{j_k}}, {y}_{n_{j_k}},{\tau}_{n_{j_k}},
{\theta}_{n_{j_k}},{\wp}_{n_{j_k}},P_{n_{j_k}})\] which converges in
$\mathscr{H}$. \hfill$\Box$

Now, we can state a theorem on the asymptotic behavior of solutions
\begin{theorem}  Let ${U}_0=({\bf{v}}^0,{ \bf{z}}^0,w^0,y^0,\tau^0,\theta^0,\wp^0,P^0)\in\mathscr{D(A)}$ and ${U}(t)$ be the solution of the
 boundary-initial-value problem (\ref{cauchy}) with $\mathcal{F}=0$.
 Then
\[ \tau(t),\ \wp(t) \to 0\ \hbox{as}\ t \to \infty\ \hbox{in}\ W_0^{1,2}(\Sigma)\quad \hbox{and}
 \quad\theta(t),\ P(t) \to 0\ \hbox{as}\ t \to \infty\ \hbox{in}\ L^2(\Sigma). \]
Moreover
\[ {\bf{v}}(t)  \to 0\ \hbox{as}\ t \to \infty\
\hbox{in}\ {\bf{W}_0^{1,2}}(\Sigma) \quad \hbox{and}
 \quad{\bf{z}}(t)    \to 0\ \hbox{as}\ t \to \infty\ \hbox{in}\ {\bf{L}}^2(\Sigma) \]\[ w(t) \to 0\ \hbox{as}\ t \to \infty\ \hbox{in}\ W_0^{1,2}(\Sigma)\quad \hbox{and}
 \quad y(t) \to 0\ \hbox{as}\ t \to \infty\ \hbox{in}\ L^2(\Sigma). \]
whenever the system
\begin{equation}
\label{subsystem}
\begin{array}{lll}
{\bf Av}+{\bf B}w+\lambda^2\bf{v}&=&0 \ \hbox{in}\ \Sigma\\  G {\bf v}+ \mathcal{H}w+\lambda^2 w&=&0 \ \hbox{in}\ \Sigma\\
  J\bf{z}&=&0\ \hbox{in}\ \Sigma\\
L\bf{z}&=&0\ \hbox{in}\ \Sigma\\
\bf{v}&=&0\ \hbox{on}\  \Gamma\\
w&=&0\ \hbox{on}\  \Gamma
\end{array}
\end{equation}
 has only the null solution.\end{theorem}
{\bf Proof.} To prove the theorem we have to study the structure of
the $\omega-$limit set. Thus, we must study the equation \be
\label{hat} \hat{\mathscr{A}}{U}=i \lambda {U}\ee for some real
number $\lambda$, where ${U}\in \mathscr{D}(\hat{\mathscr{A}})$ and
$\hat{\mathscr{A}}={\mathscr{A}}_{|\mathscr{L}}$ is the generator of
a group on $\mathscr{L}$. If ${U} \in \mathscr{L}$ then
$<\mathscr{A}{U},{U}>=0$. Under the condition \eqref{cond2}, it
follows that $\theta=P=0$ and then $\tau=\wp=0$. Thus, the
asymptotic behavior of the temperature and the chemical potential is
proved.

Now, Eq. (\ref{hat}) can be rewritten as\[ \left(
{\begin{array}{*{20}c}
  0& {\bf Id} &    0&  0&  0&  0&0&  0&   \\
      {\bf A}&  0&{\bf B}&  0 & 0&{\bf C}  &0& {\bf D} \\
0&0&0&{ Id} &    0&   0&0&  0&  \\
{ G}&  0& \mathcal{H}&   0&0  & 0&0&  0&\\
0&0&0&0 &      0& { Id}&0&  0& \\
 0&  { J}& 0&   0&\mathcal{K}  & \mathcal{L}&\mathcal{M}&  \mathcal{N}&\\
 0&0&0&0 &      0& 0&0&  Id& \\
0&  {L}& 0&   0&\mathcal{U}  & \mathcal{V}&\mathcal{W}&  \mathcal{X}&\\
\end{array}}\right)\left( {\begin{array}{*{20}c}
   {\bf v} &   \\
  {\bf z}&   \\
   w&\\
y&  \\ 0&\\
0&  \\
0&  \\
    0& \\
\end{array}}\right)=i\lambda\left( {\begin{array}{*{20}c}
   {\bf v} &   \\
  {\bf z}&   \\
   w&\\
y&  \\ 0&\\
0&  \\
0&  \\
    0& \\
\end{array}}\right).\]
Introducing the first equation into the others we obtain the system
(\ref{subsystem}).

If the system (\ref{subsystem}) has only the trivial solution, we
obtain $\omega-$limit$({U}_0) = 0$ and ${\bf{v}}(t)  \to 0\
\hbox{in}\ {\bf{W}}_0^{1,2}(\Sigma) $, $ {\bf{z}}(t)    \to 0\
\hbox{in}\ {\bf{L}}^2(\Sigma) $, $ w(t)    \to 0\ \hbox{in}\
L^2(\Sigma) $ and $ y(t)    \to 0\ \hbox{in}\ L^2(\Sigma) $ when $ t
\to \infty$.\hfill$\Box$

\section{The spatial decay estimates}
Let us consider an
unbounded body, that is we assume that $\Sigma$ is an unbounded
regular region. For the plate whose middle surface has a specific
shape, we prove that a measure associated with the thermodynamic
process decays faster than an exponential of a polynomial of second
degree.

Following \cite{Chirita0}, we consider a given
time $\mathscr{T} \in (0,\infty)$. We denote by
${\mathfrak{D}}_\mathscr{T}$ the support of the initial and boundary
data, the body force,  the heat supply and the diffusion supply on
the time interval $[0,\mathscr{T}]$, that is, the set of all ${\bf
x}\in \bar \Sigma$ such that

 $(i)$  if ${\bf x}\in \Sigma$, then \[ v^0_\alpha({\bf x})\not=0\ \ \hbox{or} \ \ z^0_\alpha({\bf x})\not=0\ \
 \hbox{or}\ \   w^0({\bf x})\not=0\ \ \hbox{or} \ \ y^0({\bf x})\not=0\ \
 \hbox{or}\ \   \tau^0({\bf x})\not=0\ \ \hbox{or} \ \ \theta^0({\bf x})\not=0\]
 \[\ \ \hbox{or} \ \  \wp^0({\bf x})\not=0\ \ \hbox{or} \ \ P^0({\bf x})\not=0\]
\be \hbox{or} \ \  f_\alpha({\bf x},t)\not=0\ \ \hbox{or} \ \ f({\bf
x},t)\not=0\ \ \hbox{or} \ \   W({\bf x},t)\not=0\ \ \hbox{or} \ \
V({\bf x},t)\not=0\ \ \ \hbox{for some }\ t\in[0,\mathscr{T}],\ee

$(ii)$  if ${\bf x} \in \p \Sigma$, we have \[ {{v}}_\alpha({\bf
x},t)\not=0\ \ \ \hbox{or} \ \ {w}({\bf x},t)\not=0\ \ \ \hbox{or} \
\ {\tau}({\bf x},t)\not=0\ \ \ \hbox{or} \ \  {\theta}({\bf
x},t)\not=0\ \ \ \hbox{or} \ \  {\wp}({\bf x},t)\not=0\ \ \
\hbox{or} \ \  {P}({\bf x},t)\not=0 \ \ \hbox{for some }\
t\in[0,\mathscr{T}].\]

Following \cite{Horgan0},  we assume that the support
${\mathfrak{D}}_\mathscr{T}$ of the initial and boundary data is
enclosed in the half-space $x_2 < 0$. We introduce the notation
$S_z$ for the open cross-section of $\Sigma$ for which $x_2 = z,\ z
\geq 0$ and whose unit normal vector is $(0, 1)$. We assume that the
unbounded set $\Sigma$ is so that $S_z$ is bounded for all finite
$z\in [0,\infty)$. We denote by $\Sigma_z$ that portion of $\Sigma$
for which $x_2
> z$.

We introduce the following function,
\be\label{func}
J(z,t)=-\int_0^t\int_{S_z}\Big( M_{2\alpha} (s)\dot v_{\alpha} (s) +2hN_{2} (s) \dot w
  (s)-\Psi_2 (s) \theta (s)-\Omega_2 (s) P (s)\Big) dx_1ds.
\ee As we want to obtain an upper estimate of the spatial decay of
the solutions, it is natural to assume that\be\label{lim} \lim_{z\to
\infty} J(z,t)=0 \ee for every finite time. Using the constitutive
equations, the divergence theorem and having in mind the definition
of $S_z$ and $\Sigma_z$ , we deduce that the function $J(z, t)$ can
be written in the form \be\label{fuct1} J(z,t)= \int_{\Sigma_z}\Big(
\mathscr{K}(t)+\mathscr{C}(t) +\mathscr{V}(t) \Big) da +
\int_0^t\int_{\Sigma_z}\mathcal{D}(s) dads,\ee where
\begin{align}\begin{split}\label{matrix}
\mathscr{K}(t)=&\half \rho\Big(I \dot v_\alpha(t)\dot v_\alpha(t)+2h \dot w^2(t)\Big), \\
\mathcal{D}(t)=& I(k_2 \theta_{,\alpha}\theta_{,\alpha}+2\hbar_2
\theta_{,\alpha}P_{,\alpha}+h_2P_{,\alpha}P_{,\alpha}) +2h
(k_2\theta^2+2\hbar_2 \theta P+h_2P^2), \\\mathscr{C}(t)=&\half I
(c\theta^2+2\kappa \theta P+ rP^2)
\end{split}\end{align} and $\mathscr{V}$ is defined by \eqref{V}.
Further, we introduce the function \be\label{en} E(z,t)=
\int_{z}^{\infty}J(r,t)dr.\ee
\begin{lemma} There exists a
positive constant  $\zeta$ depending on $k_2$ and $h_2$ such that
\be\label{vin}\int_{S_z}I(k_2\theta^{2}+h_2P^2) dx_1 \leq \zeta
\frac{\p^2 E}{\p z^2}.\ee
\end{lemma}
{\bf Proof.}   It is easy to see that \be\label{derv1} \frac{\p
E}{\p z}= -J(z,t)\ee and \be\label{derv2} \frac{\p^2 E}{\p z^2
}(z,t)= \int_{S_z}\Big( \mathscr{K}(t)+\mathscr{C}(t)
+\mathscr{V}(t) \Big) dx_1 + \int_0^t\int_{S_z}\mathcal{D}(s)
dx_1ds. \ee On the other hand, we have
\[
k_{2}\theta^{2}+h_{2}P^{2}\leq\alpha\left(\theta^{2}+P^{2}\right),
\]
with $\alpha=\max\{k_{2},h_{2}\}>0$. Moreover, being
$c\theta^{2}+2k\theta P+rP^{2}$ a positive definite quadratic form
(see \eqref{cond1}), denoting with $0<k_{\min}\leq k_{\max}$ the
smallest and the largest eigenvalues, we have
\[
k_{\min}\left(\theta^{2}+P^{2}\right)\leq c\theta^{2}+2k\theta
P+rP^{2}\leq k_{\max}\left(\theta^{2}+P^{2}\right).
\]
Then if we set $\zeta=\alpha/k_{\min}$, we get
\[
k_{2}\theta^{2}+h_{2}P^{2}\leq\zeta
k_{\min}\left(\theta^{2}+P^{2}\right)\leq\zeta\left(c\theta^{2}+2\kappa\theta
P+rP^{2}\right),
\]
so that
\[
\begin{alignedat}{1}\frac{1}{2}\int_{S_{z}}I\left(k_{2}\theta^{2}+h_{2}P^{2}\right)dx_{1} & \leq\frac{1}{2}\zeta\int_{S_{z}}I\left(c\theta^{2}+2\kappa\theta P+rP^{2}\right)dx_{1}=\zeta\int_{S_{z}}\mathscr{C}(t)dx_{1}\\
 & \leq\zeta\left(\int_{S_{z}}\left(\mathscr{K}(t)+\mathscr{C}(t)+\mathscr{V}(t)\right)dx_{1}+\int_{0}^{t}
 \int_{S_{z}}\mathcal{D}(s)dx_{1}ds\right)=\zeta\frac{\partial^{2}E}{\partial
 z^{2}}(z,t).
\end{alignedat}
\]
\hfill$\Box$

   Now, we can state the main result of this section
 \begin{theorem} There exists a
positive constant  $\xi$ depending on the constitutive constants
such that the function $E$ decays exponentially in terms of the
square of the distance $z$ from the support
${\mathfrak{D}}_\mathscr{T}$ of the external given data when $z
> \xi t $.\end{theorem}

{\bf Proof.}
 From the constitutive equations,  and the
relations \eqref{func} and \eqref{en}, we have \be\label{derv3}
\frac{\p E}{\p t}(z,t)=-
 \int_{\Sigma_z}\Big( M_{2\alpha} (t)\dot v_{\alpha} (t) +2hN_{2} (t) \dot w  (t)-\Psi_2 (t) \theta (t)-\Omega_2 (t) P (t)
 \Big) da. \ee
 It is worth remarking that

\begin{align}\begin{split}\label{tem}\int_{\Sigma_z} \Psi_2 (t) \theta (t)da=&-
 \int_{\Sigma_z}  I(k_1 \tau_{,2}+\hbar_1 \wp_{,2}+k_2
\theta_{,2}+\hbar_2 P_{,2})\theta
 da \\
=&-\int_{\Sigma_z}  I(k_1 \tau_{,2}+\hbar_1 \wp_{,2}+\hbar_2
P_{,2})\theta
 da-\half\int_{S_z}Ik_2\theta^2 dx_1
\end{split}\end{align}

and
\begin{align}\begin{split}\label{dif}\int_{\Sigma_z} \Omega_2 (t) P (t)da=&-
 \int_{\Sigma_z}  I(h_1 \wp_{,2}+\hbar_1 \tau_{,2}+\hbar_2
\theta_{,2}+h_2 P_{,2})P
 da \\
=&-\int_{\Sigma_z}  I(h_1 \wp_{,2}+\hbar_1 \tau_{,2}+\hbar_2
\theta_{,2})P
 da-\half\int_{S_z}Ih_2P^2 dx_1.
\end{split}\end{align} Thus we have

\begin{align}\begin{split}\label{derv4} \frac{\p E}{\p
t}(z,t)=&-\int_{\Sigma_z}\Big( M_{2\alpha} \dot v_{\alpha} +2hN_{2}
\dot w  +I(k_1 \tau_{,2}+\hbar_1 \wp_{,2}+\hbar_2 P_{,2})\theta
+I(h_1 \wp_{,2}+\hbar_1 \tau_{,2}+\hbar_2 \theta_{,2})P
 \Big) da
  \\
&+\half\int_{S_z}I(k_2\theta^{2}+h_2P^2) dx_1.
\end{split}\end{align}

Our next step is to estimate the time derivative of $E$ in terms of
the first two spatial derivatives of $E$. We have
$$\left|\int_{\Sigma_z}\Big( M_{2\alpha} \dot v_{\alpha} +2hN_{2}
\dot w  +I(k_1 \tau_{,2}+\hbar_1 \wp_{,2}+\hbar_2 P_{,2})\theta
+I(h_1 \wp_{,2}+\hbar_1 \tau_{,2}+\hbar_2 \theta_{,2})P
 \Big) da\right|\leq \xi J(z,t)=-\xi \frac{\p E}{\p
z}$$  where $\xi$ can be easily calculated in terms of the
constitutive constants and parameters.  This estimate together with
\eqref{vin} allow us to obtain the inequality $$\frac{\p E}{\p
t}\leq \zeta \frac{\p^2 E}{\p z^2}-\xi \frac{\p E}{\p z}.$$
Following \cite{Horgan0}, such inequality can be
further treated by the Comparison Principle \citep{Tikhonov}. On this
basis we can conclude that
$$E(z,t)\leq (\max_{s\in[0,t]}E(0,s))e^{\frac{\xi z}{2\zeta}}H(z,t)$$
where $$H(z,t)=\frac{1}{2\sqrt{\zeta \pi}}\int_0^t z
s^{-3/2}e^{-\frac{z^2}{4\zeta s}-\frac{\xi^{2}s}{4\zeta}}ds.$$ Using
the estimate discussed by \cite{Pompei} , we have
$$H(z,t)\leq \frac{2z\sqrt{\frac{\zeta t}{\pi}}}{z^2-\xi^2t^2}e^{-\frac{1}{4\zeta}(\xi^2t+\frac{z^{2}}{t})}$$
for $z>\xi t$ and $t>0$.
Thus, we find the estimate
$$E(z,t)\leq 2\max_{s\in[0,t]}E(0,s)\frac{z\sqrt{\frac{\zeta
t}{\pi}}}{z^2-\xi^2t^2}e^{-\frac{\xi^2}{4\zeta}t}e^{(\frac{\xi
z}{2\zeta}-\frac{z^{2}}{4\zeta t})}$$ for $z > \xi t$ and $t > 0$.
\hfill$\Box$

\begin{remark}Note that this result is not valid in
the case of type II theory because of $k_2=h_2 = 0$ gives $\zeta =
0$.\end{remark}

\section{Impossibility of localization in time}
In  previous sections we have proved that the  solutions of type III
theory are stable in the sense of Lyapunov, decay asymptotically and
spatially. A natural question is to ask if the decay is fast enough
to guarantee that the solution vanishes in a finite time. In fact,
when the dissipation mechanism in a system is sufficiently strong,
the localization of solutions in the time variable can hold. This
means that the decay of the solutions is sufficiently fast to
guarantee that they vanish after a finite time.

In the context of Green-Naghdi thermoelasticity of type III,
\cite{Quint00} has shown that the thermal dissipation is
not strong enough to obtain  the localization in time of the
solutions. In this section we assume the quadratic form
(\ref{positive}) positive definite and prove that the further
dissipation effects due to  diffusion is not sufficiently strong to
guarantee that the thermomechanical deformations vanish after a
finite interval of time. This means that, in absence of sources, the
only solution for the evolutive problem that vanishes after a finite
time is the null solution, that is the following theorem holds.
\begin{theorem} Let $({\bf{v}},w, \tau,\wp)$ be a solution of the system (\ref{equat}), (\ref{initial}) and  (\ref{BC})
 which vanishes after a finite time $t_0$. Then
$({\bf{v}},w,\tau,\wp)\equiv (0,0,0,0)$ for every $t \ge 0$.
\end{theorem} In order to prove this theorem, generalizing the
technique used in \cite{Quint00}, we
 show the uniqueness of solutions for the related  backward in time
 problem. Theses problems
are relevant from the mechanical point of view when we want to have
some information about what happened in the past by means of the
information that we have at this moment.

For our model, the system of equations which govern the backward in
time problem is given by

\begin{align}\begin{split}\label{back} \rho I
\ddot v_\alpha&=I(\mu \Delta
v_\alpha+(\lambda+\mu)v_{\beta,\beta\alpha}+d_1\dot
\tau_{,\alpha}+d_2\dot
\wp_{,\alpha})-2h\mu(v_\alpha+w_{,\alpha}),\cr \rho \ddot w
 &=\mu \Delta w+\mu v_{\alpha,\alpha},\cr cI\ddot \tau+\kappa I\ddot \wp&=
 I(k_1 \Delta\tau+\hbar_1 \Delta\wp-k_2 \Delta\dot \tau-\hbar_2 \Delta\dot  \wp)
 +Id_1 \dot v_{\alpha,\alpha}-2h(k_1 \tau+\hbar_1 \wp-k_2 \dot \tau-\hbar_2 \dot \wp),\cr
\kappa I\ddot \tau+r I\ddot \wp&= I(h_1 \Delta\wp+\hbar_1
\Delta\tau-\hbar_2\Delta\dot \tau-h_2\Delta \dot \wp)+Id_2 \dot
v_{\alpha,\alpha}-2h(h_1 \wp+\hbar_1 \tau-\hbar_2\dot \tau-h_2 \dot
\wp).\end{split}\end{align}

\begin{proposition} [Uniqueness]   Let   $({\bf{v}},w, \tau,\wp)$ be a solution
of the system (\ref{back}), (\ref{BC}) with null initial data and
sources. Then  $({\bf{v}},w, \tau,\wp)=(0,0,0,0)$ for every $t\ge
0$.
\end{proposition}
{\bf Proof. } Let us introduce the following functionals
\begin{align}\begin{split} \label{en1}
E_1(t)&=\half\int_\Sigma\Big(\rho I \dot {v_\alpha}\dot {v_\alpha}
+2h\rho \dot w^2+c I  \dot\tau^2+2\kappa I \dot\tau\dot \wp+rI
\dot\wp^2 +2\mathscr{V} \Big)da,\\
E_2(t)&=\half\int_\Sigma\Big(\rho I \dot {v_\alpha}\dot {v_\alpha}
+2h\rho \dot w^2-c I  \dot\tau^2-2\kappa I \dot\tau\dot \wp-rI
\dot\wp^2 +2\mathscr{U} \Big)da,\\
 E_3(t)&=\int_\Sigma I\Big(\rho {v_\alpha}\dot{v}_\alpha+2\frac{h}{I}\rho  w\dot{w}-c \dot\tau\tau-m \dot\wp\wp-\kappa(\tau\dot{\wp}
 +\dot{\tau}\wp)\\&+\frac{1}{2}
 (k_2\tau_{,\alpha}^2+2\hbar_2\tau_{,\alpha}\wp_{,\alpha}+
h_2\wp_{,\alpha}^2)+\frac{h}{I}(k_2\tau^2+ 2\hbar_2
\tau\wp+h_2\wp^2)\Big)da,\\\end{split}\end{align}
  where $\mathscr{V}$ is defined by \eqref{V} and \be \label{V2}2\mathscr{U}= I\Big(
\lambda\varepsilon_{\alpha\alpha}\varepsilon_{\beta\beta}
+2\mu\varepsilon_{\alpha\beta}\varepsilon_{\alpha\beta}-k_1\tau_{,\alpha}^2-2\hbar_1\tau_{,\alpha}\wp_{,\alpha}-
h_1\wp_{,\alpha}^2 \Big)+2h\Big(\mu
\gamma_\alpha\gamma_\alpha-k_1\tau^2- 2\hbar_1
\tau\wp-h_1\wp^2\Big).\ee

We compute now their time derivatives. By multiplying the  first
equation of (\ref{back}) by $ \dot v_\alpha$, the second one by
$\dot w$,   the third one by $\dot\tau$ and the fourth one by
$\dot\wp$, we get
\[\dot{E}_1(t)=\int_\Sigma I
\Big(k_2\dot \tau_{,\alpha}\dot \tau_{,\alpha} +2\hbar_2\dot
\tau_{,\alpha} \dot \wp_{,\alpha}+h_2\dot \wp_{,\alpha}\dot
\wp_{,\alpha}+2\frac{h}{I}(k_2\dot \tau^2 +2\hbar_2\dot \tau
P+h_2\dot \wp^2)\Big)da.\] On the other hand, if we multiply
equation of (\ref{back}) by $\dot v_\alpha$, the second one by $\dot
w$, the third one by $-\dot\tau$ and the fourth one by $-\dot\wp$,
we obtain
\[ \dot{E}_2(t)=-\int_\Sigma I \Big(2d_1  \dot \tau\dot v_{\alpha,\alpha}+2d_2  \dot \wp\dot
v_{\alpha,\alpha}+k_2\dot \tau_{,\alpha}\dot \tau_{,\alpha}
+2\hbar_2\dot \tau_{,\alpha} \dot \wp_{,\alpha}+h_2\dot
\wp_{,\alpha}\dot \wp_{,\alpha}+2\frac{h}{I}(k_2\dot \tau^2
+2\hbar_2\dot \tau P+h_2\dot \wp^2)\Big)da.\] Finally, if  we
multiply the  first equation of (\ref{back}) by $  v_\alpha$,  the
second one by $ w$, the third one by $- \tau$ and the fourth one by
$-\wp$, we get\be \dot{E}_3(t)=-\int_\Sigma \Big( I(d_1\dot \tau
v_{\alpha,\alpha}-d_1\tau_{,\alpha} \dot v_{\alpha}+d_2\dot \wp
v_{\alpha,\alpha}-d_1\wp_{,\alpha} \dot v_{\alpha}-\rho\dot
v_\alpha\dot v_\alpha - 2\frac{h}{I}\rho\dot w^{2}+c
\dot\tau^2+2\kappa\dot\tau\dot\wp+r
\dot\wp^2)+2\mathscr{U}\Big)da.\ee Moreover, a well-known identity
for type III thermoelasticity (see Eqn. (3.9) in \cite{Quint33}) for
our model becomes $$\int_\Sigma I\Big(c
\dot\tau^2+2\kappa\dot\tau\dot  \wp+m \dot\wp^2\Big)da=
  \int_\Sigma\Big(\rho I\dot {v_\alpha}\dot{v}_\alpha+2h\rho   \dot{w}^2
  -2\mathscr{U}\Big)da\,.$$
 Then we
have
\[{E}_2(t)=\int_\Sigma 2\mathscr{U}da
\]
and
\begin{eqnarray*}\dot{E}_3(t)&=&-\int_\Sigma I\Big(d_1\dot \tau
v_{\alpha,\alpha}-d_1\tau_{,\alpha} \dot v_{\alpha}+d_2\dot \wp
v_{\alpha,\alpha}-d_1\wp_{,\alpha} \dot
v_{\alpha}\Big)da\,.\end{eqnarray*} We consider the function
\begin{eqnarray*}  \mathcal{E}(t)&=&\int_0^t (\epsilon E_1(s)+ {E_2}(s)+\tilde\lambda E_3(s))ds
\\&=&\frac12 \int_0^t \int_\Sigma I \Big(\epsilon\rho  \dot {v_\alpha}\dot
{v_\alpha} +2\epsilon \frac{h}{I}\rho \dot w^2+\epsilon c
\dot\tau^2+2\epsilon\kappa  \dot\tau\dot \wp+\epsilon r \dot\wp^2
+(\epsilon+2)[\lambda\varepsilon_{\alpha\alpha}\varepsilon_{\beta\beta}
+2\mu\varepsilon_{\alpha\beta}\varepsilon_{\alpha\beta}+2\frac{h}{I}
\mu \gamma_\alpha\gamma_\alpha]\Big)da\,ds \cr&+&\frac12  \int_0^t
\int_\Sigma I \Big([\tilde\lambda k_{2}+(\epsilon-2){k}_{1}]
\tau_{,\alpha}^2 +2[\tilde\lambda\tilde{h}_{2}+(\epsilon-2)\tilde{h}_{1}]
 \tau_{,\alpha}\wp_{,\alpha}
+[\tilde\lambda h_{2}+(\epsilon-2){h}_{1}]  \wp_{,\alpha}^2\Big)da\,ds
\\&+&\frac12  \int_0^t \int_\Sigma 2h \Big([\tilde\lambda k_{2}+(\epsilon-2){k}_{1}]
\tau^2 +2[\tilde\lambda\tilde{h}_{2}+(\epsilon-2)\tilde{h}_{1}]
 \tau\wp
+[\tilde\lambda h_{2}+(\epsilon-2){h}_{1}] \wp^2\Big)da\,ds
\\&+&\tilde\lambda \int_0^t \int_\Sigma I \Big(\rho {v_\alpha}\dot{v}_\alpha+2\frac{h}{I}\rho  w\dot{w}-c \dot\tau\tau-m \dot\wp\wp-\kappa(\tau\dot{\wp}
 +\dot{\tau}\wp)\Big)da\,ds
\end{eqnarray*}
where $\epsilon$ and $\tilde\lambda$ are positive suitable constants such
that the quadratic forms

\[\int_\Sigma I \Big([\tilde\lambda k_{2}+(\epsilon-2){k}_{1}]
\tau_{,\alpha}^2 +2[\tilde\lambda\tilde{h}_{2}+(\epsilon-2)\tilde{h}_{1}]
 \tau_{,\alpha}\wp_{,\alpha}
+[\tilde\lambda h_{2}+(\epsilon-2){h}_{1}]  \wp_{,\alpha}^2\Big)da\] and
\[ \int_\Sigma I \Big([\tilde\lambda k_{2}+(\epsilon-2){k}_{1}] \tau^2
+2[\tilde\lambda\tilde{h}_{2}+(\epsilon-2)\tilde{h}_{1}]
 \tau\wp
+[\tilde\lambda h_{2}+(\epsilon-2){h}_{1}] \wp^2\Big)da
\]
are positive definite.

 By using the null initial data hypothesis and the Poincar\'e inequality we have
 \[\tilde\lambda \int_0^t \int_\Sigma I \Big(\rho {v_\alpha}\dot{v}_\alpha+2\frac{h}{I}\rho  w\dot{w}-c \dot\tau\tau-m \dot\wp\wp-\kappa(\tau\dot{\wp}
 +\dot{\tau}\wp)\Big)da\]\[\leq
\frac{\epsilon}4 \int_0^t \int_\Sigma I \Big(\rho \dot
{v_\alpha}^2+2\frac{h}{I}\rho  \dot{w}^2+ c
\dot\tau^2+2\kappa\dot\tau\dot  \wp+ m \dot\wp^2\Big)\,dv\,ds
 \]
 for any $t\leq t_0$, where $t_0$ is a positive time which depends on $\tilde\lambda$, $\epsilon$ and the constitutive coefficients.
 Therefore $\mathcal{E}(t)$ is a positive definite quadratic form for $0\leq t\leq t_0$, in particular
 \begin{eqnarray} \label{calE} \mathcal{E}(t)&\geq&\frac{1}{4} \int_0^t \int_\Sigma I \Big(\epsilon\rho  \dot {v_\alpha}\dot
{v_\alpha} +2\epsilon\frac{ h}{I}\rho \dot w^2+\epsilon c
\dot\tau^2+2\epsilon\kappa  \dot\tau\dot \wp+\epsilon r \dot\wp^2
+(\epsilon+2)[\lambda\varepsilon_{\alpha\alpha}\varepsilon_{\beta\beta}
+2\mu\varepsilon_{\alpha\beta}\varepsilon_{\alpha\beta}+2h \mu
\gamma_\alpha\gamma_\alpha]\Big)da\,ds \cr&+&\frac{1}{4}  \int_0^t
\int_\Sigma I \Big([\tilde\lambda k_{2}+(\epsilon-2){k}_{1}]
\tau_{,\alpha}^2 +2[\tilde\lambda\tilde{h}_{2}+(\epsilon-2)\tilde{h}_{1}]
 \tau_{,\alpha}\wp_{,\alpha}
+[\tilde\lambda h_{2}+(\epsilon-2){h}_{1}]  \wp_{,\alpha}^2\Big)da\,ds
\\&+&\frac{1}{4}   \int_0^t \int_\Sigma I \Big([\tilde\lambda k_{2}+(\epsilon-2){k}_{1}]
\tau^2 +2[\tilde\lambda\tilde{h}_{2}+(\epsilon-2)\tilde{h}_{1}]
 \tau\wp
+[\tilde\lambda h_{2}+(\epsilon-2){h}_{1}] \wp^2\Big)da\,ds.\nonumber
\end{eqnarray}
Moreover, recalling the null initial data assumption, we have
 \begin{eqnarray*}\dot{\mathcal{E}}(t)=(\epsilon-1)\int_0^t\int_\Sigma I \Big(k_2\dot
\tau_{,\alpha}\dot \tau_{,\alpha} +2\hbar_2\dot \tau_{,\alpha} \dot
\wp_{,\alpha}+h_2\dot \wp_{,\alpha}\dot
\wp_{,\alpha}+2\frac{h}{I}(k_2\dot \tau^2 +2\hbar_2\dot \tau
P+h_2\dot \wp^2)
\Big)dads-\int_0^t\int_\Sigma\mathscr{M}dads\,\end{eqnarray*} where

$$\mathscr{M}= I \Big(2d_1 \dot
\tau\dot v_{\alpha,\alpha}+2d_2  \dot \wp\dot
v_{\alpha,\alpha}+\tilde\lambda[d_1\dot \tau
v_{\alpha,\alpha}-d_1\tau_{,\alpha} \dot v_{\alpha}+d_2\dot \wp
v_{\alpha,\alpha}-d_1\wp_{,\alpha} \dot v_{\alpha}]\Big).$$

Choosing $0<\epsilon<1$ and using the inequality of arithmetic and
geometric means, we have
 \begin{eqnarray*}
\left|\int_0^t\int_\Sigma\mathscr{M}dads\right|&\leq&
(1-\epsilon)\int_0^t\int_\Sigma I \Big(k_2\dot \tau_{,\alpha}\dot
\tau_{,\alpha} +2\hbar_2\dot \tau_{,\alpha} \dot
\wp_{,\alpha}+h_2\dot \wp_{,\alpha}\dot
\wp_{,\alpha}+2\frac{h}{I}(k_2\dot \tau^2 +2\hbar_2\dot \tau
P+h_2\dot \wp^2) \Big)dads
\\&+&
K_1  \int_0^t \int_\Sigma (\rho I \dot {v_\alpha}\dot {v_\alpha}
+2h\rho \dot w^2)\,da\,ds +K_2 \int_0^t \int_\Sigma I(c
\dot\tau^2+2\kappa\dot\tau\dot \wp+ m
\dot\wp^2)da\,ds\\
&+&K_3   \int_0^t \int_\Sigma I
(\lambda\varepsilon_{\alpha\alpha}\varepsilon_{\beta\beta}
+2\mu\varepsilon_{\alpha\beta}\varepsilon_{\alpha\beta}+2h \mu
\gamma_\alpha\gamma_\alpha) da\,ds
\\
&+&\frac12  \int_0^t \int_\Sigma I \Big([\tilde\lambda
k_{2}+(\epsilon-2){k}_{1}] \tau_{,\alpha}^2
+2[\tilde\lambda\tilde{h}_{2}+(\epsilon-2)\tilde{h}_{1}]
 \tau_{,\alpha}\wp_{,\alpha}
+[\tilde\lambda h_{2}+(\epsilon-2){h}_{1}]  \wp_{,\alpha}^2\Big)da\,ds
\\&+&\frac12  \int_0^t \int_\Sigma 2h \Big([\tilde\lambda k_{2}+(\epsilon-2){k}_{1}]
\tau^2 +2[\tilde\lambda\tilde{h}_{2}+(\epsilon-2)\tilde{h}_{1}]
 \tau\wp
+[\tilde\lambda h_{2}+(\epsilon-2){h}_{1}] \wp^2\Big)da\,ds
\end{eqnarray*}
where the positive constants $K_i$ can be calculated by standard
methods, so that
 \begin{eqnarray}   \label{dotcalE} \dot{ \mathcal{E}}(t)&\leq&K \int_0^t \int_\Sigma I \Big(\epsilon\rho  \dot {v_\alpha}\dot
{v_\alpha} +2\epsilon \frac{h}{I}\rho \dot w^2+\epsilon c
\dot\tau^2+2\epsilon\kappa  \dot\tau\dot \wp+\epsilon r \dot\wp^2
+(\epsilon+2)[\lambda\varepsilon_{\alpha\alpha}\varepsilon_{\beta\beta}
+2\mu\varepsilon_{\alpha\beta}\varepsilon_{\alpha\beta}+2h \mu
\gamma_\alpha\gamma_\alpha]\Big)da\,ds \cr&+&K  \int_0^t \int_\Sigma
I \Big([\tilde\lambda k_{2}+(\epsilon-2){k}_{1}] \tau_{,\alpha}^2
+2[\tilde\lambda\tilde{h}_{2}+(\epsilon-2)\tilde{h}_{1}]
 \tau_{,\alpha}\wp_{,\alpha}
+[\tilde\lambda h_{2}+(\epsilon-2){h}_{1}]  \wp_{,\alpha}^2\Big)da\,ds
\\&+&K \int_0^t \int_\Sigma 2h \Big([\tilde\lambda
k_{2}+(\epsilon-2){k}_{1}] \tau^2
+2[\tilde\lambda\tilde{h}_{2}+(\epsilon-2)\tilde{h}_{1}]
 \tau\wp
+[\tilde\lambda h_{2}+(\epsilon-2){h}_{1}] \wp^2\Big)da\,ds.\nonumber
\end{eqnarray}
with $K=\max\{\frac12,K_1,K_2,K_3\}$. Inequalities (\ref{calE}) and
(\ref{dotcalE}) yield
\[ \dot{\mathcal{E}}(t)\le 4K \mathcal{E}(t)\,,\qquad 0\le t \le t_0\,.\]
This inequality and the null initial data imply
$\mathcal{E}(t)\equiv 0$ if  $0\le t \le t_0$. Reiterating this
argument on each subinterval $[(n-1)t_0 , nt_0]$ we obtain
$\mathcal{E}(t)\equiv 0$ for  $ t \geq0$.

If we take into account the definition of $\mathcal{E}(t)$,
the uniqueness result is proved. \hfill$\Box$

\section{Conclusions}
The results established in this paper can be summarized as follows:

\verb"(i)" We have derived the  linear theory of bending plate
  for thermoelastic diffusion materials of type II and III in the frame of
Green-Naghdi theory.  We have shown that the type II theory is
conservative and the solutions cannot decay with respect to time. It
is well known that, in general, the solutions of type III decay with
respect to time.

\verb"(ii)" We have proved that the problem of bending plate
  for thermoelastic diffusion materials of type  III (or type II) is well
posed. This result proves that in the motion following any
sufficiently small change in the external system, the solution of
the initial-boundary value problem is everywhere arbitrary small in
magnitude.

\verb"(iii)" We have have shown the asymptotic behaviour of
solutions of type  III model.

$(\verb"iv")$ We have have shown that the  spatial decay for the
solutions corresponding to unbounded plates, for a fixed time, at
large distance to the support ${\mathfrak{D}}_\mathscr{T}$ of the
given data is dominated by term $e^{-\frac{z}{4\zeta t}}$, where
$\zeta$ depends only on the thermal  and the diffusion coefficients
characterizing the type  III model. We can conclude that at large
distance from the support of the external given data, the spatial
decay is influenced only by the thermal and diffusion effects
arising only in type III model (see Remark 1).

$(\verb"v")$ We have derived the uniqueness of solutions of type III
model for the backward in time problem. Thus, it says the
impossibility of localization in time of the solutions. From a
thermomechanical point of view, this result says that combination of
thermal and diffusion dissipations in an  elastic bending plate is
not sufficiently strong to guarantee that the thermomechanical
deformations vanish after a finite interval of time.  Results of
this kind are a good complement to the ones we have obtained in
\cite{Passarella1,Quint00}.

We believe that by adapting the same analysis, we can  prove the
impossibility of localization of solutions in the case  of exterior
domains, even when the solutions can be unbounded, whether the
spatial variable goes to infinity. This proof is omitted for the
sake of brevity.


\section*{Acknowledgment}
Part of this work was done when the first
author visited the  Dipartimento di Ingegneria Industriale,
Universit\`a degli Studi di Salerno, in September 2015 and February
2016. He tanks their hospitality and financial support  through FARB
2014, 2015.

\end{document}